\documentclass{amsart}

\newtheorem{theorem}{Theorem}[section]
\newtheorem{corollary}[theorem]{Corollary}

\theoremstyle{definition}

\newtheorem{remark}[theorem]{Remark}
\newtheorem{example}[theorem]{Example}

\usepackage{amscd}

\newcommand{\PP}{\mathbb{P}}

\newcommand  {\shE}     {\mathcal{E}}
\newcommand  {\shF}     {\mathcal{F}}
\newcommand  {\shG}     {\mathcal{G}}

\newcommand  {\shI}     {\mathcal{I}}

\newcommand  {\shM}     {\mathcal{M}}

\newcommand  {\shL}     {\mathcal{L}}
\newcommand  {\shR}     {\mathcal{R}}

\newcommand  {\shProj}  {\mathcal{P}\!\text{\textit{roj}}\,}

\newcommand  {\aff}     {{\text{aff}}}

\newcommand  {\Fr}      {\operatorname{Fr}}

\newcommand  {\lra}     {\longrightarrow}

\renewcommand{\O}       {\mathcal{O}}

\newcommand  {\Proj}    {\operatorname{Proj}}

\newcommand  {\ra}      {\rightarrow}

\newcommand  {\Spec}    {\operatorname{Spec}}

\newcommand  {\Sym}     {\operatorname{Sym}}

\def\mydate{\number\day\space\ifcase\month \or January\or February\or March\or April\or May\or\June\or\July\or
August\or September\or October\or November\or December\fi \space\number\year}

\begin{document}

\title[Characterization of semiampleness]{A characterization of semiampleness and contractions of relative
curves}

\author[Stefan Schroeer]{Stefan Schr\"oer}
\address{Mathematisches Institut, Ruhr-Universit\"at, 
         44780 Bochum, Germany}
\curraddr{M.I.T. Department of Mathematics, 
          77 Massachusetts Avenue, Cambridge MA 02139-4307, USA}
\email{s.schroeer@ruhr-uni-bochum.de}

\subjclass{14A15, 14C20, 14H10,}


\dedicatory{Revised version}

\begin{abstract}
I give a cohomological  characterization
of semiample line bundles. The result  is a 
generalization of both the Fujita--Zariski
Theorem on semiampleness  and the Grothendieck--Serre
Criterion for ampleness.
As an application of the Fujita--Zariski Theorem I
characterize  contractible curves
in 1-dimen\-sional families.
\end{abstract}

\maketitle

\section*{Introduction}

The Fujita--Zariski Theorem  asserts that
a line bundle 
$\shL$ that  is  ample on its base locus   is \emph{semiample}.
Semiampleness  means that  a multiple 
$\shL^{\otimes n}$, $n>0$ is globally generated. For discrete base locus the result 
goes back to Zariski
(\cite{Zariski 1962}, Thm.~6.2), and the general form is due to Fujita
(\cite{Fujita 1983}, Thm.~1.10).
This note contains two applications of the Fujita--Zariski Theorem.

The first section contains  a  generalization 
of both the Fujita--Zariski Theorem  and the  
cohomological criterion for ample\-ness due to Gro\-then\-dieck--Serre.
The  result is the following characterization:
A line bundle $\shL$ is semiample if and only if the modules
$H^1(X,\shI\otimes\Sym\shL)$ are finitely generated  over the ring 
$\Gamma(X,\Sym\shL)$ for every coherent 
ideal $\shI\subset\O_B$. Here $B\subset X$ is the stable base locus of $\shL$.
This gives a positive answer to Fujita's question  
(\cite{Fujita 1983},  1.16) whether it
is possible to weaken the assumption in the
Fujita--Zariski Theorem.

In the second section I generalize  results of 
Piene \cite{Piene 1974} and Emsalem \cite{Emsalem 1974}.
They used the Fujita--Zariski Theorem to obtain sufficient conditions
for contractions in  normal arithmetic surfaces. Our result is  a 
characterization of  
contractible curves in 1-dimensional families
over  local noetherian rings in terms of complementary closed subsets.
This also sheds some light on  the
noncontractible curve constructed  by Bosch, L\"utkebohmert, and Raynaud 
(\cite{Bosch; Luetkebohmert; Raynaud 1990},  chap.~6.7). 
For proper  normal algebraic surfaces,
similar results appear in 
\cite{Schroeer 2000}.

\section{Characterization of semiampleness}
Throughout this section, 
$R$ is a noetherian ring,
$X$ is  a proper 
$R$-scheme, and  
$\shL$ is an invertible 
$\O_{X}$-module. According to the   Grothendieck--Serre Criterion
(\cite{Grothendieck 1961a}, Prop.~2.6.1)
$\shL$ is ample if and only if for each coherent 
$\O_{X}$-module 
$\shF$ there is an integer 
$n_0>0$ so that 
$H^1(X,\shF\otimes\shL^{\otimes n})=0$  for all 
$n>n_0$. 
Let me reformulate this in terms of graded modules. For a coherent
$\O_X$-module $\shF$, set 
$$
H_*^p(\shF,\shL)=
H^p(X,\shF\otimes\Sym\shL)=
\bigoplus_{n\geq 0} H^p(X,\shF\otimes \shL^{\otimes n}).
$$ 
This is a graded 
module over the graded ring 
$\Gamma_*(\shL)=\Gamma(X,\Sym \shL)$. 
The Gro\-then\-dieck--Serre Criterion takes the form: 
$\shL$ is ample if and only if the 
modules 
$H^1_*(\shF,\shL)$ are finitely generated  over the ring
$\Gamma_0(\shL)=\Gamma(\O_{X})$ for all coherent 
$\O_{X}$-modules 
$\shF$.
In this form  it generalizes to the semiample case.
Following Fujita
\cite{Fujita 1983}, we define the \emph{stable base locus} $B\subset X$ of 
$\shL$ to be the intersection of the base loci of 
$\shL^{\otimes n}$ for all
$n>0$. We regard it as a closed subscheme with reduced scheme structure.

\begin{theorem}
\label{semiample}
Let 
$B\subset X$ be the stable base locus of 
$\shL$. Then the following are equivalent:
\renewcommand{\labelenumi}{(\roman{enumi})}
\begin{enumerate}
\item 
The invertible sheaf $\shL$ is semiample.
\item
The 
modules
$H^p_*(\shF,\shL)$ are finitely generated
over the ring $\Gamma_*(\shL)$   for each coherent 
$\O_{X}$-module 
$\shF$ and all integers 
$p\geq 0$.
\item
The 
modules 
$H^1_*(\shI,\shL)$ are finitely generated over the ring
$\Gamma_*(\shL)$  for each coherent 
ideal $\shI\subset \O_B$.
\end{enumerate}
\end{theorem}

\proof 
The implication (i)$\Rightarrow$(ii) is well known, and 
(ii)$\Rightarrow$(iii) is trivial.
To prove 
(iii)$\Rightarrow$(i) we assume that 
$\shL$ is not semiample. According to the Fujita--Zariski 
Theorem the restriction 
$\shL_B$ is not ample. By the Grothendieck--Serre Criterion 
there is a coherent ideal
$\shI\subset \O_B$ with 
$H^1(X,\shI\otimes\shL^{\otimes n})\neq 0$   for infinitely many 
$n>0$. Thus 
$H^1_*(\shI,\shL)$ is not finitely generated  over 
$\Gamma_0(\shL)$. Since 
$B\subset X$ is the  stable base locus, the maps 
$\Gamma(X,\shL^{\otimes n})\ra\Gamma(B,\shL_B^{\otimes n})$ vanish for all 
$n>0$. Consequently, the irrelevant ideal 
$\Gamma_+(\shL)\subset\Gamma_*(\shL)$ annihilates 
$H^1_*(\shI,\shL)$, which is therefore  
not finitely generated  over 
$\Gamma_*(\shL)$.
\qed

\medskip
Sommese \cite{Sommese 1978} introduced a 
quantitative version of semiampleness: Let $k\geq 0$ be an integer;
a semiample invertible sheaf
$\shL$  is called
\emph{$k$-ample} if  the fibers of the canonical morphism 
$f:X\ra\Proj\Gamma_*(\shL)$ have dimension 
$\leq k$. For example, 
$0$-ampleness  means ampleness.

\begin{theorem}
\label{k-ample}
Let 
$\shL$ be a semiample invertible 
$\O_{X}$-module.  Then
$\shL$ is 
$k$-ample if and only if  the modules
$H^{k+1}_*(\shF,\shL)$ are finitely generated  over the ground
ring $R$  for all coherent 
$\O_{X}$-modules 
$\shF$.
\end{theorem}

\proof
Set 
$Y=\Proj \Gamma_*(\shL)$ and let 
$f:X\ra Y$ be the corresponding contraction.
Suppose 
$\shL$ is 
$k$-ample. Choose 
$n_0>0$ so that $\shL^{\otimes n_0}=f^*(\shM)$ for some ample invertible 
$\O_Y$-module 
$\shM$. Put 
$\shG=
\shF\otimes(\shL\oplus\shL^{\otimes 2}\oplus\ldots\oplus\shL^{\otimes n_0})$.
Choose 
$m_0>0$ with 
$H^p(Y,R^qf_*(\shG)\otimes\shM^{\otimes m})=0$ for 
$p>0$, $q\leq k+1$,  and 
$m>m_0$. Consequently, the edge map  
$H^{k+1}(X,\shG\otimes \shL^{\otimes mn_0}) \ra H^0(Y,R^{k+1}f_*(\shG)\otimes\shM ^{\otimes m})$
in the spectral sequence 
$$
H^p(Y,R^qf_*(\shG)\otimes\shM ^{\otimes m})\Longrightarrow 
H^{p+q}(X,\shG\otimes \shL^{\otimes mn_0})
$$
is injective for 
$m>m_0$. The fibers of 
$f:X\ra Y$ are at most 
$k$-dimensional, so 
$R^{k+1}f_*(\shG)=0$. Thus $H^{k+1}(X,\shF\otimes\shL^{n})=0$  for all
$n>n_0m_0$.
 
Conversely, assume that the condition holds.
Seeking a contradiction we suppose that some 
fiber of 
$f:X\ra Y$  has dimension 
$>k$. Using
\cite{Kleiman 1967} we   find a coherent 
$\O_{X}$-module 
$\shF$ with 
$R^{k+1}f_*(\shF) \neq 0$. Replacing 
$\shL$ by a suitable multiple, we have 
$\shL=f^*(\shM)$ for some ample invertible 
$\O_Y$-module $\shM$. Passing to a higher multiple if necessary, 
$H^p(Y,R^qf_*(\shF)\otimes\shM^{\otimes n})=0$ holds for  
$p>0$, $q\leq k$,  and 
$n>0$. 
Then the edge map
$H^{k+1}_*(X,\shF\otimes\shL^{\otimes n})\ra H^0_*(Y,R^{k+1}f_*(\shF)\otimes\shM^{\otimes n})  $ is
surjective for $n>0$. Choose a global section 
$s\in\Gamma(Y,\shM^{\otimes n})$  for some  $n>0$ so that the open subset
$Y_s\subset Y$ contains the set of associated points for 
$R^{k+1}f_*(\shF)$. Then  
$s\in\Gamma_*(\shM)$  is  not a zero divisor for
$H^0_*(R^{k+1}f_*(\shF) ,\shM)$. 
It follows that 
$H^0_*(R^{k+1}f_*(\shF) ,\shM)$ is nonzero for infinitely many degrees.
Consequently, the same holds for 
$H^{k+1}_*(\shF,\shL)$, which is therefore not finitely generated over 
$R$.
\qed

\begin{remark}
For a \textit{vector bundle} 
$\shE$, it might happen that 
$\O_{\PP(\shE)}(1)$ is semiample, whereas
$\Sym^n(\shE)$ fails to be globally generated for all 
$n>0$. For example, let 
$k$ be an algebraically closed field of characteristic 
$p>0$, and 
$X$ be a   smooth proper curve of genus 
$g>p-1$ so that the absolute Frobenius
$\Fr_X:H^1(\O_{X})\ra H^1(\O_{X})$ is zero.
For an example see 
\cite{Hartshorne 1977}, p.~348, ex.~2.14. 
Let  
$D\subset X$ be a divisor of degree 
$1$. According to the commutative diagram 
$$
\begin{CD}
H^0(\O_X)     @>>> H^0(\O_D) @>>> H^1(\O_{X}(-D)) @>>> H^1(\O_{X})    \\
@V\Fr_X^* VV     @V\Fr_X^* VV @V\Fr_X^* VV @VV\Fr_X^*=0 V\\
H^0(\O_X)     @>>> H^0(\O_{pD}) @>>> H^1(\O_{X}(-pD)) @>>> H^1(\O_{X}),  
\end{CD}
$$
the 
$p$-linear map 
$\Fr_X^*:H^1(\O_{X}(-D))\ra H^1(\O_{X}(-pD))$ is not injective. 
Hence there is a nontrivial extension 
$$
0 \lra \O_{X}  \lra \shE  \lra \O_X(D)  \lra 0
$$
whose Frobenius pull back 
$\Fr_X^*(\shE)$ splits. The surjection 
$\shE\ra\O_X(D)$ gives a section 
$A\subset \PP(\shE)$ representing 
$\O_{\PP(\shE)}(1)$ with 
$A^2=1$
(\cite{Hartshorne 1977}, Prop.\ 2.6, p.\ 371). 
The Fujita--Zariski Theorem implies that 
$\O_{\PP(\shE)}(1)$ is semiample, and we obtain a birational contraction 
$\PP(\shE)\ra Y$.  It is easy to see that  the exceptional set 
is an integral curve 
$R\subset \PP(\shE)$ which has degree 
$p$ on the ruling. Hence 
$\PP(\shE)\ra Y$ does  not restrict to closed embeddings on the fibers of 
$\PP(\shE)\ra X$. Consequently, 
$\Sym^n(\shE)$ is not globally generated at any point 
$x\in X$.
\end{remark}

\section{Contractions of relative curves}

Throughout this section, 
$R$ is a local noetherian ring, and 
$X$ is  a proper 
$R$-scheme with 1-dimensional closed fiber 
$X_0\subset X$. Then  all fibers of the structure morphism
$X\ra\Spec(R)$ are at most 1-dimensional. For example, $X$ could
be a flat family of curves.

A \emph{Stein factor} of 
$X$ is a proper 
$R$-scheme 
$Y$ together with a proper morphism 
$f:X\ra Y$ so that 
$\O_Y\ra f_*(\O_X)$ is bijective (compare 
\cite{Kleiman 1966}, sec.~5). Our objective is to describe
the set of all Stein factors for a given 
$X$.

Let 
$C_i$, 
$i\in I$ be the finite collection of all 1-dimensional integral
components of the closed fiber
$X_0$. A subset 
$J\subset I$ yields a subcurve 
$C=\bigcup_{i\in J}C_i$. We call such a curve
$C\subset X$ \emph{contractible} if there is a Stein factor 
$f:X\ra Y$ so that 
$f(C_i)$ is a closed point  if and only if
$i\in J$. 
According to 
\cite{Grothendieck 1961a}, Theorem 5.4.1, 
a Stein factor is determined up to isomorphism by its restriction
$f_0:X_0\ra Y_0$. The task now is to determine the contractible curves 
$C\subset X$.
It follows from  \cite{Piene 1974} and 
\cite{Emsalem 1974} that all curves 
$C\subset X$ are contractible provided that the ground ring
$R$ is henselian. In particular this holds if 
$R$ is complete. On the other hand,
a noncontractible curve is discussed in 
\cite{Bosch; Luetkebohmert; Raynaud 1990}, chapter 6.7.

We seek to describe contractible curves $C\subset X$ in terms of  
complementary closed subsets 
$D\subset X$. We need a definition:
Suppose 
$D\subset X$ is a  closed subset of codimension 
$\leq 1$.
Let 
$R\subset R^\wedge$ be the completion with respect to the maximal ideal, 
$X'$ the normalization of $X\otimes_R R^\wedge$, and
$C_i',C',D'\subset X'$ the preimages of  $C_i,C,D\subset X$, respectively.  Let
$h:X'\ra Z'$ be the contraction of all
$C_i'\subset X_0'$ disjoint from  
$C' $. We call   
$D$  \emph{persistent} if 
$h(D') \subset Z'$  has codimension   
$\leq 1$.

\begin{example}
Suppose 
$R$ is a discrete valuation ring with residue field 
$k$ and fraction field 
$K$. Let 
$X$ be the proper $R$-scheme obtained from 
$X'=\PP^1_R$ by identifying the closed points 
$0,\infty\in\PP^1_k$. Then the closure 
$D\subset X$ of the point
$0\in\PP^1_K$ is not persistent.
\end{example}

\begin{theorem}
\label{contractible}
Suppose 
$J\subset I$ is  a subset so that the curve 
$C=\bigcup_{i\in J}C_i$ is  connected.
Then $C\subset X_0$ is contractible if and only if there
is a persistent closed subset 
$D\subset X$ of codimension 
$\leq 1$ disjoint from 
$C$ and intersecting each irreducible component
$C_i\subset X_0$ with 
$i\not\in J$.
\end{theorem}

\proof
Assume  that $C$ is contractible. The corresponding  contraction 
$f:X\ra Y$ maps 
$C$ to a single point. Let 
$V\subset Y $ be an affine open neighborhood of 
$f(C)$. Set $U=f^{-1}(V)$ and 
$D=X- U$. Clearly 
$D\cap C=\emptyset$. 
Furthermore, 
$D\cap C_i\neq \emptyset$ for 
$i\not\in J$; otherwise $f(C_i)$ would be
a proper curve contained in the affine scheme
$V$, which is absurd.
Let 
$X',Y'$ be the normalizations of 
$X\otimes_R R^\wedge, Y\otimes_R R^\wedge$, respectively. The induced morphism 
$f':X'\ra Y'$ is the contraction  
 of the preimage 
$C'\subset X'$ of
$C$. The preimage 
$V'\subset Y'$ of 
$V$ is affine, so $Y- V$ is of codimension 
$\leq 1$ (\cite{SGA 6} II, 2.2.6). Hence the preimage
$D'\subset X'$ of 
$D$ is of codimension 
$\leq 1$. Obviously, the same holds if we contract the preimages 
$C_i'\subset X'$ of 
$C_i$ disjoint from 
$C'$. Thus 
$D\subset X$ is of codimension
$\leq 1$ and persistent.

Conversely, assume the existence of such a subset 
$D\subset X$. Set 
$U=X- D$. We claim that the affine hull 
$U^\aff=\Spec\Gamma(U,\O_{X})$ is of finite type over 
$R$ and that the canonical morphism 
$U\ra U^\aff$ is proper.

Suppose this for a moment. Then 
$U\ra U^\aff$ contracts 
$C$ and is a local isomorphism near each
$x\in U_0- C$.
Choose for each 
$x\in X_0- C$ an affine open neighborhood 
$U_x\subset X$ of 
$x$ disjoint to the exceptional set of 
$U\ra U^\aff$. Then 
$U_x\cap U\ra U^\aff$ is an open embedding. It is easy to see that
the schemes 
$U_x\bigcup_{U_x\cap U} U^\aff$,  $x\in X_0- C$ and 
$U^\aff$ form an open cover of a proper 
$R$-scheme 
$Y$. The induced morphism 
$f:X\ra Y$ is the desired contraction.

It remains to verify the claim. Let $R\subset R^\wedge$ be the completion.
According to \cite{SGA 1}, VIII Corollary 3.4, the scheme $U^\aff$ is 
of finite type if and only if $U^\aff\otimes_R R^\wedge$ is of finite type.
Furthermore, $U\ra U^\aff$ is proper if and only if if  is proper
after tensoring with $R^\wedge$ (\cite{SGA 1}, VIII Cor.~4.8).
Since $U^\aff\otimes_R R^\wedge=(U\otimes_R R^\wedge)^\aff$ 
by \cite{Grothendieck 1967}, Proposition 21.12.2,
it suffices to prove the claim under the additional assumption 
that $R$ is complete.

Now each curve in 
$X_0$ is contractible.
Observe that the contraction of
$C$ does not change 
$U^\aff$, so we  can as well  assume that 
$C$ is empty.
Now our goal is to prove that 
$U$ is affine.
Since 
$R$ is complete, hence universally japanese, the normalization 
$X'\ra X$ is finite. Using Chevalley's Theorem
(\cite{Grothendieck 1961}, Thm.~6.7.1),
we reduce the problem to the case that 
$X$ is normal. 
Now the irreducible components of $X$ are the connected components.
Treating them separately we may assume that
$X$ is connected. Contracting the curves 
$C_i$ contained in  
$D$ we can assume that 
$D_0$ is finite and intersects each 
$C_i$.
If 
$D= X$  or $D=\emptyset$ there is nothing to prove.
Assume that 
$D\subset X$ is of codimension 1, in other words a Weil divisor.
The problem is that it might not  be  Cartier. 
To overcome this, consider the graded quasicoherent 
$\O_{X}$-algebra 
$\shR=\bigoplus_{n\geq 0}\O_{X}(nD)$. The  graded subalgebra
$\shR'\subset \shR$  generated by 
$\shR_1=\O_{X}(D)$   is of finite type over 
$\O_{X}$. Set 
$X'=\shProj(\shR')$ and let 
$g:X'\ra X$ be the structure morphism. Then 
$g$ is projective and 
$\O_{X'}(1)$ is a 
$g$-very ample invertible 
$\O_{X'}$-module.
The canonical maps 
$D:\O_{X}(nD)\ra\O_{X}((n+1)D)$ induce a homomorphism 
$\shR'\ra\shR'$ of degree one, hence  a section 
$s:\O_{X'}\ra\O_{X'}(1)$.  It  follows from the definition of 
homogeneous spectra that  
$s$ is bijective over 
$U$ and vanishes on 
$g^{-1}(D)$. Thus the corresponding Cartier divisor 
$D'\subset X'$ representing 
$\O_{X'}(1)$ has support 
$g^{-1}(D)$.

Let 
$A\subset X_0'$ be a closed integral subscheme of dimension 
$n>0$. If 
$g(A)\subset X_0$ is a curve, then 
$A$ is not contained in 
$D'$ but  intersects 
$D'$. Hence 
$D'\cdot A>0$. 
If 
$g(A)\subset X$ is a point, then 
$\O_{A}(1)$ is ample, so 
$(D')^n\cdot A>0$. By the Nakai criterion for ampleness we conclude that 
$\O_{X'}(1)$ is ample on its base locus. Now the  
Fujita--Zariski Theorem tells us that 
$\O_{X'}(1)$ is semiample. It follows that 
$U\simeq X'- D'$ is  affine. This finishes the proof.
\qed

\medskip
Let us consider the special case that the total space
$X$ is a normal surface. Replacing $R$ by $\Gamma(X,\O_X)$, we are in the
following situation: Either 
$R$ is   a discrete valuation ring, such  that 
$X\ra\Spec(R)$ is a flat deformation of  $X_0$.
Or 
$R$ is a local normal 2-dimensional ring, hence 
$X\ra\Spec(R)$ is the birational contraction of 
$X_0$. In either case  we call a  Weil divisor 
$H\in Z^1(X)$   \emph{horizontal} if it is a sum of 
prime divisors not supported by
$X_0$. 

Suppose 
$J\subset I$ is a subset with 
$C=\bigcup_{i\in J}C_i$ connected.
Let 
$V\subset X_0$ be the union of all 
$C_i$ disjoint from
$C$.

\begin{corollary}
\label{surfaces}
Notation as above. Then 
$C\subset X_0$ is contractible if and only if there is a horizontal 
Weil divisor 
$H\subset X$ disjoint from 
$C$ with the following property: For each 
$C_i$, 
$i\not\in J$, either 
$H$ intersects $C_i$, or 
$H$ intersects a connected component   
$V'\subset V$ with
$V'\cap C_i\neq\emptyset$.
\end{corollary}

\proof
Suppose 
$C\subset X_0$ is contractible. Let 
$D\subset X$ be a persistent Weil divisor as in Theorem 
\ref{contractible}. Then its horizontal part 
$H\subset D$ satisfies the above conditions. 
Conversely, assume there is a horizontal Weil divisor 
$H\subset X$ as above. It follows that 
$D=H+V$ is a persistent Weil divisor disjoint from 
$C$  intersecting each 
$C_i$ with $i\not\in J$. Thus 
$C\subset X_0$ is contractible.
\qed


\end{document}